\documentclass[11pt]{article}

\usepackage{epsfig}
\usepackage{graphicx}
\usepackage[latin1]{inputenc}

\usepackage{amsfonts,amsthm,bbm,amssymb,epsf,amsmath,%
latexsym,amscd,amsfonts,enumerate,amsthm,supertabular,color,url}

\usepackage{subfigure,psfrag}

\newcommand{\p}{\partial}

\usepackage{algorithm}
\usepackage{algorithmic}

\newcommand{\mpr}[1]{\relax}
\newcommand{\mprup}[1]{\relax}
\newcommand{\mprdown}[1]{\relax}

\newtheorem{theorem}{Theorem}[section]

\newtheorem{definition}[theorem]{Definition}

\newtheorem{remark}[theorem]{Remark}

\newtheorem{testproblem}{Test Problem}

\newcommand{\pf}{{\noindent \bf Proof.} }
\newcommand{\eop}{{ \vrule height7pt width7pt depth0pt}\par\bigskip}

\def\sp{:}
\def\ie{{\sl i.e.},\ }

\def\tri{\triangle}

\newcommand{\RR}{\mathbb{R}}

\newcommand{\PP}{\mathbb{P}}

\newcommand{\sps}{S^{1,2}_{5,0}(\tri)}
\newcommand{\spsh}{S^{1,2}_{5,0}(\tri^h)}

\newcommand{\wG}{{\widetilde G}}
\newcommand{\hu}{{\hat u}}

\DeclareMathOperator{\supp}{supp}

\DeclareMathOperator{\st}{star}

\begin{document}

\title{$C^1$ Quintic Splines 
on Domains Enclosed by Piecewise
Conics and Numerical Solution of Fully Nonlinear Elliptic Equations}%

\author{Oleg Davydov\thanks{Department of
Mathematics, University of Giessen,
    Department of Mathematics,
    Arndtstrasse 2, 
    35392 Giessen,
    Germany, {\tt oleg.davydov@math.uni-giessen.de}}
\and Abid Saeed\thanks{Department of Mathematics,
Kohat University of Science and Technology,
Kohat, Pakistan, {\tt abidsaeed@kust.edu.pk}}}

\maketitle

\begin{abstract}
We introduce bivariate $C^1$ piecewise quintic finite element spaces  for 
curved domains enclosed by piecewise conics satisfying homogeneous boundary conditions,
construct local bases for them using Bernstein-B\'ezier techniques, and demonstrate 
the effectiveness of these finite elements for the numerical solution of 
the Monge-Amp\`ere equation over curved domains by B\"ohmer's method.
\end{abstract}

\section{Introduction\label{intro}}\mprup{intro}

Piecewise polynomials on curved domains bounded by piecewise algebraic curves and surfaces is a
promising but little studied tool for data fitting and solution of partial differential
equations. 
Since  implicit algebraic surfaces are a well-established modeling technique
in CAD \cite{ImplSurf}, we are interested in developing isogeometric schemes
\cite{HCB05} for domains with such boundaries,  where the geometric models of the boundary 
are used exactly in the
form they exist in a CAD system rather than undergoing a remeshing to fit into the traditional
isoparametric finite element approach. 

In this paper we continue the work started in \cite{DKS}, where $C^0$ splines vanishing on a
piecewise conic boundary have been introduced.
In contrast to both the isoparametric curved finite elements and the isogeometric analysis of 
\cite{HCB05}, our approach 
does not require parametric patching on curved subtriangles, and therefore does not depend on the invertibility
of the Jacobian matrices of the nonlinear geometry mappings. Therefore our finite elements remain piecewise polynomial
everywhere in the physical domain. 

This approach allows to incorporate conditions of higher smoothness in Bernstein-B\'ezier
form standard for the theory and practice of smooth piecewise polynomials on polyhedral domains
\cite{LSbook}. It turns out however that imposing boundary conditions make the
otherwise well understood spaces of e.g.\ bivariate $C^1$ macro-elements on triangulations 
significantly more complex. Even in the simplest case of a polygonal domain, the
dimension of the space of splines vanishing on the boundary is dependent on its geometry, with 
consequences for  the construction of stable bases (or stable minimal determining sets) 
\cite{DSa12,DSa13}.

In this paper we suggest a local basis defined through a minimal determining set for the space of $C^1$
piecewise quintic polynomials vanishing on a piecewise conic boundary and apply the resulting
finite element space to the numerical solution of the fully nonlinear Monge-Amp\`ere equation on
 domains with such boundary. The latter is done within the framework of Böhmer's method 
 \cite{Boehmer08} which we applied previously on polygonal domains \cite{DSa13}. The results are
 based in part on the thesis of the second named author \cite{Abid_thesis}.
 
It is important to mention that  the  isoparametric 
approach to $C^0$ curved elements is problematic
 when finite element spaces of $C^1$ or higher smoothness are sought, see the remarks in 
 \cite[Section~4.7]{BrennerScott}. A successful $C^1$ quintic construction of this type developed in
\cite{MB93} seems difficult to extend to higher smoothness or higher polynomial degree.

Remarkably, the standard Bernstein-B\'ezier  
techniques for dealing with piecewise polynomials on triangulations 
\cite{LSbook,Sch15} as well as recent
optimal assembly algorithms \cite{AAD11,AAD15,ADS} for high order elements 
are carried over to the spaces used here without significant loss of
efficiency, see \cite{DKS}.

The paper is organized as follows. The spaces of $C^1$ piecewise polynomials on domains with
piecewise conic boundary are introduced in Section~2, whereas Section~3 presents our
construction of a local basis for the main space of interest $\sps$. Section~4 briefly
summarizes Böhmer's method for fully nonlinear elliptic equations and presents a number of 
numerical experiments for the Monge-Amp\`ere equation on smooth domains, including a circular
domain, an elliptic domain, and piecewise conic domains with $C^1$ and $C^2$ boundaries.

\section{$C^1$ piecewise polynomials on  piecewise conic domains\label{spacesC0}}\mprup{spaces}

We first recall from \cite{DKS} the assumptions on a domain $\Omega$ and its 
triangulation $\tri$ with curved pie-shaped triangles at the boundary. 

Let $\Omega\subset\RR^2$ be a bounded curvilinear polygonal domain with 
$\Gamma=\partial \Omega=\bigcup_{j=1}^n\overline{\Gamma}_j$, 
where each $\Gamma_j$ is an open arc of an algebraic curve of
at most second order
(\ie either a straight line or a conic). For simplicity we assume that 
$\Omega$ is simply connected.
Let $Z=\{z_1,\ldots,z_n\}$ be the set of the endpoints of all
arcs numbered counter-clockwise such that $z_j,z_{j+1}$ are the
endpoints of $\Gamma_j$, $j=1,\ldots,n$, with $z_{j+n}=z_j$. 
Furthermore, for each $j$ we denote by $\omega_j$ the internal 
angle between the tangents $\tau^+_j$ and $\tau^-_j$ to $\Gamma_j$
and $\Gamma_{j-1}$, respectively, at $z_j$. We assume that $\omega_j>0$ for all $j$.

\begin{figure}[htbp!]
\centering
\begin{tabular}{c}
\includegraphics[height=0.3\textwidth]{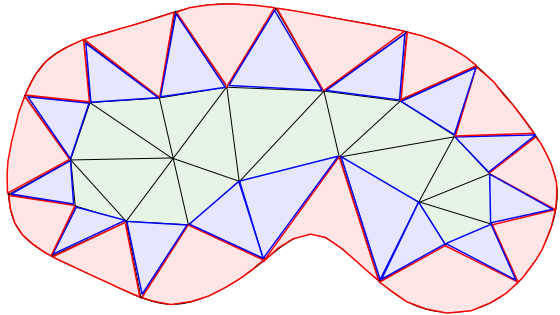} 
\end{tabular}
\caption{
A triangulation of a curved domain with ordinary triangles (green), 
pie-shaped triangles (pink) and buffer triangles (blue).
}\label{curved_tri}
\end{figure}

Let $\tri$ be a \emph{triangulation} of $\Omega$, \ie a 
subdivision of $\Omega$ into triangles, where 
each triangle $T\in\tri$ has at most one edge replaced with a curved segment of 
the boundary $\partial\Omega$, and the intersection of any pair of the triangles
is either a common vertex or a common (straight) edge if it is non-empty.
The triangles with a curved edge are said to be \emph{pie-shaped}.
Any triangle $T\in\tri$ that shares at least one edge with a pie-shaped triangle
is called a \emph{buffer} triangle, and the remaining triangles are
\emph{ordinary}. We denote by $\tri_0$, $\tri_B$ and $\tri_P$ the sets of all
ordinary, buffer and pie-shaped triangles of $\tri$, respectively, such that
$\tri=\tri_0\cup\tri_{B} \cup \tri_P$
is a disjoint union, see Figure~\ref{curved_tri}.
Let $V,E,V_I,E_I,V_B,E_B$ denote the set of all vertices, 
all edges, interior vertices, interior edges, boundary vertices and boundary  edges,
respectively.

For each  $j=1,\ldots,n$, let $q_j\in\mathbb{P}_2$ be a polynomial such that 
$\Gamma_j\subset\{x\in\RR^2\sp q_j(x)=0\}$, where $\mathbb{P}_d$ denotes the space of all bivariate polynomials of total degree
at most $d$. By changing the sign of $q_j$ if needed, we ensure that
 $\p_{\nu_x} q_j(x)<0$ for all $x$ in the interior of
$\Gamma_j$, where $\nu_x$ denotes the unit outer normal to the boundary at $x$,
and $\p_{a}:=a\cdot\nabla$ is the directional derivative with respect to a vector $a$. Hence,
$q_j(x)$ is positive for points in $\Omega$ near the boundary segment $\Gamma_j$.
We assume that $q_j\in \mathbb{P}_1$ if $\Gamma_j$ is a straight interval. Clearly,
$q_j$ is an irreducible quadratic polynomial if $\Gamma_j$ is a  genuine conic arc and in all
cases
\begin{equation}\label{qgrad}
\nabla q_j(x)\ne 0\quad\text{if}\quad x\in \Gamma_j. 
\end{equation}

Following \cite{DKS} we assume that $\tri$ satisfies the following conditions:
\begin{itemize}
\item[(a)] $Z=\{z_1,\ldots,z_n\}\subset V_B$.
\item[(b)] No interior edge has both endpoints on the boundary.
\item[(c)] No pair of pie-shaped triangles shares an edge.
\item[(d)] Every $T\in\tri_P$ is star-shaped with respect to its interior vertex $v$.
\item[(e)]  For any  $T\in\tri_P$ with its curved side on $\Gamma_j$,
$q_j(z)>0$ for all $z\in T\setminus \Gamma_j$.
\end{itemize}
It can be easily seen that (b) and (c) are achievable by a slight modification of
a given triangulation, while (d) and (e) hold for sufficiently fine triangulations.

For any $d\ge1$ we set
\begin{align*}
S^1_d(\tri)&:=\{s\in C^1(\Omega)\sp s|_T\in
\mathbb{P}_{d+i},\;T\in\tri_i,\;i=0,1\},
\quad \tri_1:=\tri_P\cup\tri_B,
				 \\
S^{1,2}_{d,I}(\tri)&:=\{s\in  S^1_d(\tri)\sp s \text{ is twice differentiable at any }v\in V_I\}
				 ,\\
S^{1,2}_{d,0}(\tri)&:=\{s\in S^{1,2}_{d,I}(\tri)\sp s|_\Gamma = 0\}. %
\end{align*}

As in \cite{DKS} we use Bernstein-B\'ezier techniques  to obtain a local
basis for $\sps$ with the help of a minimal determining set. 

Recall (see \cite{LSbook}) that the 
bivariate Bernstein polynomials with respect to a non-degenerate triangle $T=\langle v_{1},v_{2},v_{3}\rangle$
with vertices $v_{1},v_{2},v_{3}\in\mathbb{R}^2$ are defined by 
 \begin{equation*}
B_{ijk}^{d}(v):= \frac{d!}{i!j!k!}b_{1}^{i}b_{2}^{j}b_{3}^{k},\quad i+j+k=d,  
\end{equation*}
where $b_{1},b_{2},b_{3}$ are the barycentric coordinates of $v$, that is the unique coefficients of the
expansion $v=\sum _{i=1}^{3}b_{i}v_{i}$ with $\sum _{i=1}^{3}b_{i}=1$. The Bernstein polynomials
form a basis for $\mathbb{P}_d$, and the coefficients  $c_{ijk}$ in the \emph{BB-form}
expansion 
\begin{equation}\label{BBform}
p=\sum _{i+j+k=d}c_{ijk}B_{ijk}^{d},\quad p\in \mathbb{P}_d,
\end{equation}
are called the \emph{BB-coefficients} of $p$.
They are conveniently indexed by the elements of the set
\begin{equation}\label{domainpoints}
D_{d,T}:=\left\lbrace \xi _{ijk} =\frac{iv_{1}+jv_{2}+kv_{3}}{d}:\;i+j+k=d,\;i,j,k \geq 0 \right\rbrace
\end{equation}
of so-called  \emph{domain points}, such that $B_{\xi}^{d}:=B_{ijk}^{d}$ and $c_{\xi}:=c_{ijk}$ when
$\xi=\xi_{ijk}\in D_{d,T}$.
We will also use the notation $D^{d,T}_2(v)$ for the subset of $D_{d,T}$ consisting
of the six domain points closest to a vertex  $v$ of $T$, in particular 
$$
D^{d,T}_2(v_1)=\{\xi_{d,0,0},\xi_{d-1,1,0},\xi_{d-1,0,1},\xi_{d-2,2,0},
\xi_{d-2,0,2},\xi_{d-2,1,1}\}.$$

The continuity and $C^1$-smoothness of piecewise polynomials are expressed as follows. 
Given two triangles $T=\langle v_{1},v_{2},v_{3}\rangle$ and 
$\tilde T=\langle v_{4},v_{3},v_{2}\rangle $ sharing an edge 
$e=\langle v_{2},v_{3}\rangle$, let $p$ and $\tilde p$
be two polynomials of degree $d$ written in the BB-form
$$
p=\sum _{i+j+k=d}c_{ijk}B_{ijk}^{d} \quad\text{ and }\quad
\tilde p=\sum _{r+s+t=d}\tilde{c}_{rst}\tilde{B}_{rst}^{d},$$
where $B_{ijk}^{d}$ and $\tilde B_{rst}^{d}$ are the Bernstein polynomials with respect to $T$ and
$\tilde T$, respectively. Then $p$ and $\tilde p$ join  continuously along $e$ 
if and only if their 
BB-coefficients over $e$ coincide, i.e. 
\begin{equation}\label{smoothnessc0}
\tilde{c}_{0jk}=c_{0kj},\quad\text{for all } j+k=d.
\end{equation}
Moreover, the condition for $C^1$ smoothness across $e$ is that \eqref{smoothnessc0} holds along with 
\begin{equation}\label{smoothnessc1}
\tilde{c}_{1jk}=b_1c_{1,k,j}+b_2c_{0,k+1,j}+b_3c_{0,k,j+1},\quad j+k=d-1,
\end{equation}
where  $(b_{1},b_{2},b_{3})$ are the barycentric coordinates of $v_4$ relative to $T$.

A finite set $\Lambda$ of linear functionals $\lambda:\sps\to\RR$ is said to
be a \emph{determining set} if 
$$
\lambda(s)=0 \quad\forall  \lambda \in \Lambda\quad 
\Longrightarrow \quad s=0,$$
and $\Lambda$ is a \emph{minimal determining set (MDS)} if there is no smaller determining
set. In other words, a determining set is a spanning set of the dual space
$(\sps)^*$, and an MDS is a basis of $(\sps)^*$. Any 
MDS $\Lambda$ uniquely determines a \emph{basis} 
$\{s_\lambda:\lambda\in\Lambda\}$  of $S$ by duality, such that
$\lambda(s_\mu)=\delta_{\lambda,\mu}$, for all $\lambda,\mu\in\Lambda$,
and any spline $s\in S$ can be uniquely written in the form
$s=\sum_{\lambda\in\Lambda}c_\lambda s_\lambda$, with
$c_\lambda=\lambda(s)\in\mathbb{R}$.

To explain what we mean by a local basis we need some further definitions, 
compare \cite{D01,DKS}.
The \emph{$\ell$-star} of a set $A\subset\Omega$ with respect to $\tri$
is given by
$$
\st^1(A)=\st(A):=\bigcup\{T\in\tri: T\cap A\ne\emptyset\},
\;\st^\ell(A):=\st(\st^{\ell-1}(A)),\; \ell\ge2.$$
A set $\omega\subset\Omega$ is said to be a
\emph{supporting set} of a linear functional $\lambda\in (\sps)^*$ if 
$\lambda(s)=0$ for all $s\in \sps$ such that $s|_\omega=0$.
Given an MDS $\Lambda$, we define for each $T\in\tri$ the set
$\Lambda_T:=\{\lambda\in\Lambda: T\subset\supp s_\lambda\}$,
where $\{s_\lambda:\lambda\in\Lambda\}$ is the basis of $\sps$ dual to $\Lambda$.
Thus, $\lambda\in\Lambda_T$ if and only if for a spline $s\in S$, 
$s|_T$ depends on the coefficient $c_\lambda=\lambda(s)$. The \emph{covering number}
$\kappa_\Lambda$ of an MDS $\Lambda$
is the maximum number of elements in  $\Lambda_T$ for all $T\in\tri$.

\begin{definition}\label{localMDS}
\rm
A minimal determining set $\Lambda$ for $\sps$ is said to be
\emph{$\ell$-local} if there is a family of supporting sets $\omega_\lambda$ of 
$\lambda\in\Lambda$ such that $\omega_\lambda\subset \st^{\ell}(T)$
for any $T\in\tri$ such that $\lambda\in\Lambda_T$. If $\Lambda$ is 
\emph{$\ell$-local} for some $\ell$, then the dual basis 
$\{s_\lambda :\lambda\in\Lambda\}$ is said to be \emph{local}.
\end{definition}

It is easy to check, see \cite[Lemma 4.3]{DKS}, that  if $\Lambda$ is 
$\ell$-local, then the basis functions $s_\lambda$ are locally supported in the sense that 
$\supp s_\lambda\subset\st^{2\ell+1}(T)$ for some triangle $T\in\tri$.

\section{A local basis for $\sps$}\label{basis}

In this section we describe a minimal determining set $\Lambda$ for $\sps$, 
which in turn defines a basis $\{s_\lambda :\lambda\in\Lambda\}$
as explained in the previous section. For the sake of simplicity we describe the basis 
under the following additional assumption:
\begin{itemize}
\item[(f)] All boundary edges are curved.
\item[(g)] No pair of buffer triangles shares an edge.
\end{itemize}
In fact we have implemented our bases also for the case where some boundary edges are straight.
(It is used in Test Problem~\ref{C1domain} in Section~4.3.)
In  this case we nevertheless assume that the triangle attached to a straight boundary edge is
ordinary, and no pie-shaped triangle shares an edge with it, as in Figure~\ref{InitMeshC1}. 
 A description of this construction would take too much space 
because it has to include the handling of the boundary vertices on ordinary polygonal 
triangulations along the lines of \cite{DSa12,DSa13}, and so we avoid this by assuming (f).
Similarly, allowing buffer triangles to share edges, or equivalently, allowing more than one
buffer triangle attached to a single boundary vertex would produce additional degrees of
freedom on and near these edges and around the boundary vertex, also requiring the techniques
of  \cite{DSa12,DSa13}.

We denote by $V_B^1$ the set of those 
boundary vertices $v\in V_B$ where the boundary $\partial\Omega$ has a well-defined
tangent, that is either $v\notin Z$, or $\omega_j= \pi$ if $v=z_j$ for some $j=1,\ldots,n$.
In addition, $E_{P,B}$ denotes the set of all edges shared by  a pie-shaped and a
buffer triangle. We also set $E_I^0:= E_I\setminus E_{P,B}$.

 Since splines in $\sps$ are 
polynomials of degree $d=5$  on the triangles $T\in\tri_0$, we can write 
these polynomials in BB-form \eqref{BBform}, 
\begin{equation}\label{BB0}
s|_T=\sum _{\xi\in D_{5,T}}c_{\xi}B_{\xi}^{5},\quad s\in\sps.
\end{equation}
and define for each $\xi=\xi_{ijk}\in D_{5,T}$ a functional $\lambda_\xi\in (\sps)^*$
that picks the BB-coefficient $c_{ijk}$ in \eqref{BBform}. With the usual convention
(see \cite{LSbook}) we identify the functional $\lambda_\xi$ with the domain point $\xi$
and speak of an MDS as a set $M\subset\overline{\Omega}$. Thanks to
\eqref{smoothnessc0} for domain points $\xi$ at vertices or on the edges of the
subtriangulation $\tri_0$ it does not matter which triangle in $\tri_0$ containing $\xi$ is used to evaluate the
BB-form of a spline $s\in \sps$. The union $D_{5,\tri_0}=\cup_{T\in\tri_0}D_{5,T}$ 
forms the standard set of domain points (and corresponding functionals $\lambda_\xi$) 
associated with $\tri_0$. Following the standard construction of an MDS for the space
$S^{1,2}_5(\tri_0)$ with only ordinary triangles \cite{LSbook}, we define the following subsets
of $D_{5,\tri_0}$. For each $v\in V_I$ we choose a triangle 
$T_v= \langle v_1,v_2,v_3\rangle\in\tri_0$ attached to $v$, such that $v_1=v$, and set 
$M_v:=D^{5,T_v}_2(v)=\{\xi_{500},\xi_{410},\xi_{401},\xi_{320},\xi_{302},\xi_{311}\}\subset D_{5,T_v}$.
For each edge $e \in E_I^0$, let 
$T_e : = \left\langle v_1,v_2,v_3\right\rangle$ 
be a triangle in $\tri_0$ attached to the edge 
$e= \left\langle v_2,v_3 \right\rangle $ and let $M_e:=\{\xi _{122}\}\subset D_{5,T_e}$. 
Clearly, $\omega_\xi:=T_v$ (resp.\ $\omega_\xi:=T_e$) is a supporting set for any functional 
$\lambda_\xi$ with $\xi\in M_v$ (resp.\ $\xi\in M_e$).

\begin{figure}[htbp!]
\centering
\begin{tabular}{c}
\psfrag{p}{$\xi$}
\includegraphics[height=0.40\textwidth]{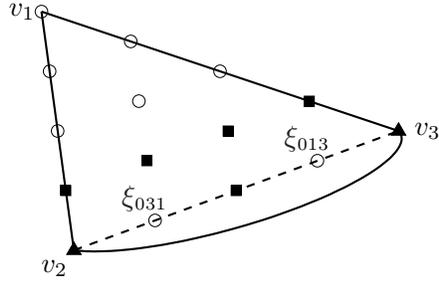} 
\end{tabular}
\caption{The set $D^\ast_{4,T}$ for a pie-shaped triangle $T$ and
domain points in $M_T^P$ (black squares), and $M_{v_2}^P\cup M_{v_3}^P$ (black triangles)
under the assumption that $v_2, v_3 \in V_B^1$ and $T=T_{v_2}=T_{v_3}$.
}\label{PMDS}
\end{figure}

For each $T \in \tri_P$, with its curved edge $e$ given by the equation $q(x)=0$, 
where $q\in \mathbb{P}_2 \backslash \mathbb{P}_1$ is irreducible and normalized so that
$q(v)=1$ for the interior vertex $v$ of $T$, we notice that 
by B\'ezout theorem
$$
\{ s\in \mathbb{P}_{6}\sp s|_{e}=0\}
=q\mathbb{P}_{4} :=\{qp\sp p\in \mathbb{P}_{4}\}.$$
Let $T^\ast$ denote the triangle obtained by joining the boundary vertices of $T$
by a straight line segment (see the dashed line in Figure~\ref{PMDS}).
Since the Bernstein polynomials $B_{ijk}^{4}$, $i+j+k=4$, w.r.t. $T^\ast$ form a basis for 
$\mathbb{P}_{4}$ it is obvious that the set 
$$
\left\lbrace qB_{ijk}^{4}\sp i+j+k=4\right\rbrace $$ 
is a basis for $q\mathbb{P}_{4}$. The set of domain points 
of degree $4$ over $T^\ast$ will be denoted $D^\ast_{4,T}$. Even though
the set $D^\ast_{4,T}$ formally coincides with $D_{4,T^\ast}$, the linear 
functionals associated with the domain points are different.
Namely, each $\xi \in D^\ast_{4,T}$ represents a linear functional 
$\lambda_\xi$ on $\sps$
which picks the  coefficient $c_\xi$ in the expansion 
\begin{equation}\label{BBpie}
s|_T=q\sum_{\xi \in D^\ast_{4,T}}c_\xi B_\xi ^{4},\qquad 
s\in \sps.
\end{equation}
Assuming that $v_1,v_2,v_3$ are the vertices of a pie-shaped triangle $T\in\tri_P$,  
with $v_1 \in V_I$, we set 
$M_T^{P}:=\{\xi _{130}, \xi _{121},\xi _{112},\xi _{103},\xi _{022}\}\subset D^\ast_{4,T}$,
see Figure~\ref{PMDS} where the points in $M_T^{P}$ are marked as black squares.
Clearly, $\omega_\xi:=T$ is a supporting set for $\lambda_\xi$.
The vertices $v_2,v_3$ of $T$ are shared by a pair of pie-shaped triangles and may belong to $V_B^1$.
For each  $v\in V_B^1$ let $M_v^P:=\{v\}\subset D^\ast_{4,T_v}$, where $T_v$ is one of the two pie-shaped 
triangle attached to $v$, and the corresponding functional is $\lambda_v$ that picks the
respective coefficient $c_v$ in \eqref{BBpie} for $T=T_v$. A supporting set for $\lambda_v$ is
given by $\omega_v:=T_v$.

\begin{figure}[htbp!]
\centering
\begin{tabular}{c}
\includegraphics[height=0.40\textwidth]{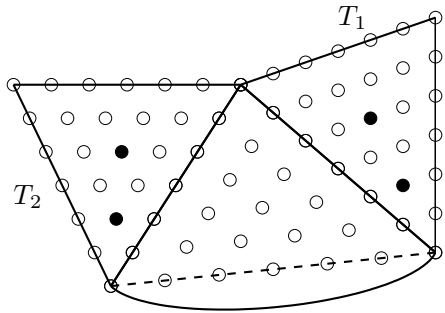} 
\end{tabular}
\caption{
The domain points in the sets $M_{T_1}^{B}, M_{T_2}^{B}$ for the buffer triangles $T_1,T_2$
are marked with black dots.
}\label{BMDS}
\end{figure}

For each $T  = \left\langle v_1,v_2,v_3\right\rangle $ in $\tri_{B}$, where $v_1 \in V_B$, 
let $M_T^{B}:=\{\xi_{411}, \xi_{222}\}\subset D_{6,T}$, see Figure~\ref{BMDS}.
As usual, the functional $\lambda_\xi$ identified with $\xi\in D_{6,T}$
picks the coefficient $c_\xi$ in the BB-form expansion of $s|_T\in\PP_6$,
\begin{equation}\label{BBB}
s|_T=\sum _{\xi\in D_{6,T}}c_{\xi}B_{\xi}^{6},\quad s\in\sps,
\end{equation}
and $\omega_\xi:=T$ is a supporting set for any $\xi\in M_T^{B}$.

\begin{remark}\label{remark1ch5}\rm
Let $T:=\left\langle v_1,v_2,v_3\right\rangle \in \tri _P$ with $v_1 \in V_I$.
Then $s|_{T}=qp \in \PP_6$ for some $p \in P_4$, where the equation $q(x)=0$ 
represents the curved edge of $T$, with an irreducible quadratic polynomial $q$ such that
$q(v_1)=1$. We can write all three polynomials $s|_{T},q,p$ in BB-form  
with respect to $T^\ast$,
\begin{equation}\label{qBB}
q = q_{110}B_{110} ^2 + q_{101}B_{101} ^2 + q_{011}B_{011} ^2 + B_{200} ^2
\end{equation}
(where we used the fact that $q(v_2)=q(v_3)=0$),
$$
s\vert _{T}= \sum _{i+j+k=6} a_{ijk}B_{ijk} ^6,\quad 
p = \sum _{i+j+k=4} c_{ijk}B_{ijk} ^4.$$
If the coefficients $c_{ijk}$ are known, then $a_{ijk}$ can be computed by multiplying the
expansions for $p$ and $q$, see the explicit formulas in \cite[Eq.~(35)]{DKS}, where a
different numeration of the vertices of $T$ is used. Moreover,
the coefficients $c_{ijk}$ can be obtained from $a_{ijk}$ in a stable way  
\cite[Lemma 4.6]{DKS}. To compute $c_{ijk}$ we may write down the identity
\begin{equation}\label{BBprodPie}
\Big(  \sum _{i+j+k=4} c_{ijk}B_{ijk} ^4\Big) 
\Big(  \sum _{i+j+k=2} q_{ijk}B_{ijk} ^2\Big)= \sum _{i+j+k=6} a_{ijk}B_{ijk} ^6 
\end{equation}
as a linear system with respect to the vector of unknown coefficients $c_{ijk}$,
$i+j+k=4$. It is easy to check that the matrix of this system has a block structure, and by
singling out the six rows of the system corresponding to the domain points in 
$D^{6,T^*}_2(v_1)$  we obtain a non-singular triangular linear system for the coefficients $c_{ijk}$
corresponding to the domain points in $D^{4,T^*}_2(v_1)$, namely 
$$ \left[ \begin{array}{cccccc}
 1 & 0 & 0 & 0 & 0 & 0\\ 
\frac{1}{3}q_{110} & \frac{8}{15} & 0  & 0 & 0 & 0 \\
\frac{1}{3}q_{101} & 0 & \frac{8}{15} & 0  & 0 & 0 \\
0 & \frac{1}{5}q_{110} & 0 & \frac{2}{5}  & 0 & 0 \\
\frac{1}{15}q_{011} & \frac{2}{15}q_{101} & \frac{2}{15}q_{110} & 0  & \frac{4}{15} & 0 \\
0 & 0 & \frac{1}{5}q_{101} & 0  & 0 & \frac{2}{5} 
\end{array}
\right] 
\cdot
\left[ \begin{array}{c}
c_{400}\\ c_{310}\\ c_{301}\\ c_{220}\\ c_{211}\\ c_{202}\end{array}\right]
=\left[ \begin{array}{c}
a_{600}\\ a_{510}\\ a_{501}\\ a_{420}\\ a_{411}\\ a_{402} 
\end{array}
\right] .$$  
Thus, we can compute the BB-coefficients $\{c_\xi\sp \xi \in D^{4,T^*}_2(v_1)\}$ of $p$
by using only the BB-coefficients $\{a_\xi\sp \xi \in D^{6,T^*}_2(v_1)\}$ of $s|_T$.
\end{remark}

\begin{theorem}\label{MDSCh5}
The set 
\begin{equation}
M := \bigcup _{v\in V_I} M_v \cup \bigcup_{e \in E_I^0} M_e 
\cup \bigcup _{v \in V_B^1} M_v^P\cup \bigcup _{T\in \tri_{P}} M _T^P 
\cup \bigcup _{T\in \tri_{B}} M _T^{B} 
\end{equation}
is a 1-local minimal determining set for the space $\sps$.
\end{theorem}

\pf
Following the standard scheme \cite{LSbook} we assign some arbitrary values $c_\xi\in\RR$ to 
$\lambda_\xi(s)$, for all $\xi\in M$, and show that all other coefficients $c_\xi$ of $s \in \sps$ on all
triangles $T\in\tri$ in the form \eqref{BB0}, \eqref{BBpie} or \eqref{BBB} depending on 
the type of $T$, can be determined from them consistently. The success of this process will
show that $M$ is an MDS. In the same time we will keep track how far the influence of a
coefficient  $c_\xi$ for any $\xi\in M$ extends, to check the locality of this MDS.

It is easy to see that the set 
$$
M_0:=\bigcup _{v\in V_I} M_v \cup \bigcup_{e \in E_I^0} M_e $$
is a 1-local MDS for the Argyris finite element space 
$$
S^{1,2}_{d}(\tri_0):=\{s\in  S^1_d(\tri_0)\sp
 s \text{ is twice differentiable at any vertex }v \text{ of }\tri_0\}$$
as shown in \cite[Theorem 6.1]{LSbook}. 

Let $v \in V_I$ be shared by some pie-shaped triangle $T \in \tri_P$. Then there are also two
buffer triangles $T_1,T_2\in\tri_B$ attached to $v$, see Figures~\ref{curved_tri} and \ref{BMDS}. We know that 
$M_v=D^{5,T_v}_2(v)\subset D_{5,T_v}$ for some $T_v\in\tri_0$. By \cite[Lemma 5.10]{LSbook}
and the degree raising formulas of \cite[Theorem 2.39]{LSbook},
$M_v$ consistently determines the BB-coefficients of $s|_{T\cup T_1\cup T_2}$ in 
$D^{6,T^*}_2(v)\cup D^{6,T_1}_2(v)\cup D^{6,T_2}_2(v)$. 
For the pie-shaped triangle $T$ we need to go one more step and compute 
the BB-coefficients in $D^{4,T^*}_2(v)$ of the polynomial $p\in\PP_4$ such that $s|_T=pq$,
where the equation $q(x)=0$ represents the curved edge of $T$. This can be done uniquely 
by solving the triangular linear system described in Remark~\ref{remark1ch5}. 

Let $e=\langle v_2,v_3\rangle\in E^0_I$ be shared by an ordinary triangle 
$T_e:=\left\langle v_1,v_2,v_3\right\rangle \in \tri _0$ and a buffer triangle 
$T=\left\langle v_4,v_3,v_2\right\rangle \in \tri _{B}$. Assuming that the BB-coefficients of
$s|_{T_e}$ for all domain points in $D_{5,T_e}$ have been computed as described above, we 
can use degree raising to write  $s|_{T_e}$ as a polynomial of degree six, and
obtain its BB-coefficients for all domain points in $D_{6,T_e}$. By using the continuity and
$C^1$ smoothness conditions \eqref{smoothnessc0}, \eqref{smoothnessc1} we can then compute the
BB-coefficients of $s|_{T}$ for all those domain points $\xi_{ijk}$ in $D_{6,T}$, for which 
$i\in\{0,1\}$. Some of them have already been computed at the previous step, namely those that belong to 
$D^{6,T}_2(v_2)\cup D^{6,T}_2(v_3)$. It is known that no inconsistencies arise this way, see for
example the proof of \cite[Theorem 6.1]{LSbook}. We thus obtain three new BB-coefficients of 
$s|_{T}$ corresponding to the domain points $\xi_{033},\xi_{132},\xi_{123}\in D_{6,T}$.

Let $v \in V_B$ and let $T_1,T_2\in\tri_P$ be the two pie-shaped triangles attached to $v$,
with the curved edges given by $q_1(x)=0$ and $q_2(x)=0$, respectively. Let $p_1,p_2\in\PP_4$
be such that $s|_{T_i}=p_iq_i$, $i=1,2$. Since $s$ is continuously differentiable at $v$
and $q_1(v)=q_2(v)=0$, we have $\nabla s(v)=p_1(v)\nabla q_1(v)=p_2(v)\nabla q_2(v)$. 
If $v \in V_B\setminus V_B^1$, then the vectors $\nabla q_1(v)$ and $\nabla q_2(v)$ are linearly
independent, and it follows that $p_1(v)=p_2(v)=0$, that is $c_v=0$ in \eqref{BBpie} for both
$T_1$ and $T_2$.
We now assume that $v \in V_B^1$. Then $\nabla q_1(v)=\alpha \nabla q_2(v)$ for some real
$\alpha\ne 0$, which implies $p_2(v)=\alpha p_1(v)$.
Let $T_1=T_v$ be the triangle in the definition of $M^P_v$, 
in particular the functional $\lambda_v$ is evaluated as $\lambda_v(s)=p_1(v)$. Thus, the value
$c_v$ in \eqref{BBpie} for $T=T_1$ is known because  $M^P_v$ is part of the MDS $M$, and the
value of the BB-coefficient of $p_2$ at the same point $v$ is $\alpha c_v$. To compute
$\alpha$, we just need to compare the components of the vectors $\nabla q_1(v)$ and 
$\nabla q_2(v)$, which is easy to do by using the BB-forms \eqref{qBB} of $q_1,q_2$ with respect to 
$T_1^*,T_2^*$, respectively. 
 
Let $T_1=\left\langle v_1,v_2,v_3\right\rangle \in \tri _{B}$ with $v_1 \in V_B$ and 
$e=\left\langle v_1,v_3\right\rangle \in E_{P,B}$, and let
$T_2:=\left\langle v_3,v_4,v_1\right\rangle \in \tri _{P}$ share the edge $e$ with $T$ and has
its curved edge defined by the equation $q(x)=0$.  Let us write the polynomials 
$s|_{T_1}$, $s|_{T_2}$ and $p\in\PP_4$ in $s|_{T_2}=pq$ in the BB-form as
$$
s|_{T_1}=\sum _{\xi \in D_{6,T_1}} \tilde c_{\xi}B_{\xi} ^6,\quad
s|_{T_2}=\sum _{\xi \in D^\ast_{6,T_2}} a_{\xi}B_{\xi} ^6,\quad
p=\sum_{\xi \in D^\ast_{4,T_2}}c_\xi B_\xi ^{4}.$$
Since the domain point $\xi_{103}\in D^\ast_{4,T_2}$ belongs to $M_{T_2}^P$ and the
coefficients $c_\xi$ for all other $\xi\in D^\ast_{4,T_2}\cap e$ have been determined above, 
$s|_e$ is completely determined, and the BB-coefficients $a_{\xi}$ for all
$\xi\in D_{6,T_2}\cap e$ can be found by the multiplication of $p|_e$ by $q|_e$. Hence the 
smoothness conditions \eqref{smoothnessc0} and \eqref{smoothnessc1} across $e$ 
give us in particular the equation 
\begin{align*}
a_{114}&=b_1\tilde c_{501}+b_2\tilde c_{411}+b_3\tilde c_{402}\\
&=b_1a_{105}+b_2\tilde c_{411}+b_3a_{204}, 
\end{align*}
where $\left(b_1,b_2,b_3 \right) $ are the barycentric coordinates of $v_4$ w.r.t.\ $T_1$,
which determines $a_{114}$ since $\xi_{411}\in D_{6,T_1}$ belongs to $M_{T_1}^B$.
Moreover, comparing the   coefficients of $B_{114}^6$ on both sides of 
  \eqref{BBprodPie} leads to the equation
$$%
 15a_{114}=q_{110}c_{004}+4q_{101}c_{013}+4q_{011}c_{103}, 
$$%
and hence $c_{013}$ can be computed from the already known BB-coefficients as
$$
c_{013}=\tfrac{1}{4q_{101}}\big(15a_{114}
-q_{110}c_{004}-4q_{011}c_{103}\big).$$
Note that $q_{101}\ne0$ thanks to \eqref{qgrad}.  Similarly, $c_{031}$ is computed using 
the same argumentation involving the buffer triangle attached to $v_4$. This completes the
computation of the BB-form of $p$. By multiplying it with $q$ we get the missing coefficients
of the BB-form of $s|_{T_2}$, and by the smoothness conditions across $e$ the BB-coefficients
$\tilde c_{312}$ and $\tilde c_{213}$ of $s|_{T_1}$. The remaining unset BB-coefficients of 
$s|_{T_1}$ are obtained in the same way by using the pie-shaped triangle sharing the edge
$\langle v_1,v_2\rangle$ with $T_1$. 

A close inspection of the above argumentation shows that $M$ is 1-local in the sense 
of Definition \ref{localMDS}.
\eop

\begin{figure}[htbp!]
\centering
\begin{tabular}{c}
\includegraphics[height=0.55\textwidth]{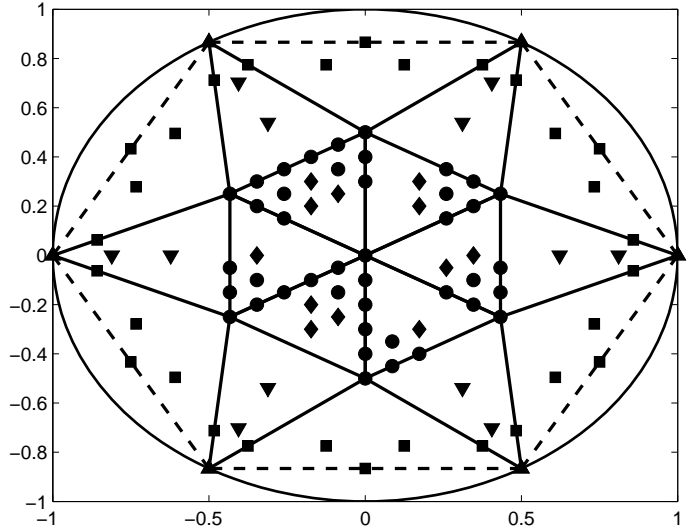} 
\end{tabular}
\caption{
Example of the  MDS of Theorem~\ref{MDSCh5} for the space $\sps$ over  a triangulation of 
a circular domain $\Omega$.
}\label{MDSC1ch5}
\end{figure}

An example of the  MDS of Theorem~\ref{MDSCh5} for the space $\sps$ over a triangulation of a
circular disk is depicted in Figure~\ref{MDSC1ch5}, 
where the points in the sets $\bigcup _{v \in V_I}M_v$, 
$\bigcup _{e \in E_{I}^0}M_e$, 
$\bigcup _{v \in V_B^1}M^P_v$, 
$\bigcup _{T \in \tri _P}M_T^P$  and 
$\bigcup _{T \in \tri _{B}}M_T^{B}$ are 
marked as black dots, diamonds, triangles,  squares and downward pointing triangles,
 respectively. Note that $V_B^1=V_B$ in this example.

\section{Numerical solution of fully nonlinear elliptic equations}
\label{numex}

To evaluate the performance of our construction 
of $C^1$ elements for curved domains we implemented B\"ohmer's method 
for fully nonlinear equations using $\sps$ as the finite element approximation space. 

\subsection{B\"ohmer's method}

We consider the Dirichlet problem,  
\begin{align}
\label{e01}
\text{find } u:\Omega\to\RR \text{ such that }G(u)= 0 \text{ and }
 u|_{\partial\Omega} = \phi , 
\end{align}
for a second order differential operator of the form
$G(u)=\wG(\cdot,u,\nabla u,\nabla ^{2}u)$, 
where  $\wG=\wG(w)$, $w=(x,z,p,r)\in\mathbb{R} \times\mathbb{R} \times 
\mathbb{R}^{2}\times \mathbb{R}^{2\times 2}$ is a real valued function
defined on a domain $\widetilde\Omega \times \Gamma$ such that 
$\overline{\Omega}\subset \widetilde\Omega\subset\mathbb{R}^2$ and 
$\Gamma \subset \mathbb{R} \times \mathbb{R}^{2}\times \mathbb{R}^{2\times 2}$,
where $\nabla u,\nabla ^{2}u$ denote the gradient and the Hessian of $u$, respectively.
The operator $G$ is said to be \emph{elliptic} in a subset 
$\widetilde\Gamma \subset \widetilde\Omega \times \Gamma$ if 
the matrix $[ \frac{\partial \wG}{\partial r_{ij}}(w)]_{i,j=1}^2 $ 
is well defined and positive definite for all $w\in \widetilde\Gamma$ \cite{BoehmerBook,GTbook}.  
Under certain assumptions, including the exterior sphere condition
for $\partial\Omega$, the continuity of 
$\phi:\partial\Omega\to\RR$ and sufficient smoothness of $\wG$,
the problem \eqref{e01} has a unique solution $u \in C^{2}(\Omega) \cap C(\overline{\Omega})$
if $\widetilde\Gamma = \widetilde\Omega \times \Gamma$ \cite[Theorem 17.17]{GTbook}. 
 
The most famous example of a fully nonlinear elliptic operator which is neither quasilinear nor semilinear 
\cite[p.~80]{BoehmerBook} is the \emph{Monge-Amp\`ere operator} 
$G(u):=\det(\nabla^{2}u)-g$, where $g:\Omega\to\RR$ satisfies $g(x)>0$ for all $x\in\Omega$. In this case 
$\widetilde\Gamma=\widetilde\Omega\times\mathbb{R}\times
\mathbb{R}^2\times\{r\in \mathbb{R}^{2\times 2}:r\hbox{ is positive definite}\}$.
Under the assumptions that $\partial\Omega$ is $C^3$ and $g\in C^2(\overline{\Omega})$ 
there exists a unique \emph{convex} solution $u$ of \eqref{e01} such that
$u \in C^{2,\alpha}(\overline{\Omega})$ for all $\alpha<1$ \cite[Theorem 17.22]{GTbook}. 
References to further results about the existence and uniqueness of the solution of \eqref{e01} 
can be found in \cite[Section 2.5.7]{BoehmerBook}.

Many fully nonlinear elliptic operators and corresponding equations $G(u)=0$ are 
important for applications, see \cite{BoehmerBook}. 
Several numerical methods have been proposed in the literature, 
in particular finite difference \cite{DG06,OA08} and finite element type methods 
\cite{GA,Boehmer08,NBGS,DSa13,XF09-1,LP13,Neilan14}.
To the best of our knowledge  however, no method has been tested 
before on non-polygonal domains.

Finite element spaces $S^h_0\subset C^1(\overline{\Omega})$ satisfying homogenous boundary 
conditions on $\Omega$, where $h$ is the maximum diameter of the underlying  partition $\Delta^h$, 
can be employed in 
\emph{Böhmer's method} \cite{Boehmer08,BoehmerBook} for the problem \eqref{e01}.
For a fixed $h>0$, let $u^h_0:\Omega\to\RR$ be an initial guess satisfying the boundary condition
$u^h_0|_{\partial\Omega} = \phi$. We generate a sequence of functions 
$\{u^h_k\}_{k\in \mathbb{N}}$ by the Newton type method
\begin{equation}\label{newton}
u^h_{k+1}=u^h_k-u^h,\quad k=0,1,\ldots,
\end{equation}
where $u^h\in S^h_0$ is the Galerkin approximation of the linear elliptic problem
\begin{equation}\label{linear}
G'(u^h_k)u=G(u^h_k),
\end{equation}
that is $u^h\in S^h_0$ is determined by the equations
\begin{equation}\label{fem}
\hbox{$(G'(u^h_k)u^h,v^h)_{L^2(\Omega)}=(G(u^h_k),v^h)_{L^2(\Omega)} 
\quad \forall v^{h}\in S^h_0$},
\end{equation}
where $( \cdot,\cdot)_{L^2(\Omega)}$ denotes the usual inner product in $L^2(\Omega)$,
and $G'(u^h_k)$ is the linearization of the operator $G$ at $u^h_k$ given  by
\begin{equation}\label{frechet}
G'(u^h_k)u=\frac{\partial \wG}{\partial z}(w^h_k)u+ \sum _{i=1}^{2} \frac{\partial \wG}{\partial
p_{i}}(w^h_k)\frac{\partial u}{\partial x_i}
+\sum _{i,j=1}^{2}\frac{\partial \wG}{\partial r_{ij}}(w^h_k)\frac{\partial^2 u}{\partial x_ix_j},
\end{equation}
with $w^h_k(x):=(x,u^h_k(x),\nabla u^h_k(x),\nabla^2u^h_k(x))$, $x\in\Omega$. 
Clearly, \eqref{fem} can be reformulated into the standard  weak form of the Galerkin method: 
Find $u^h \in S^h_0$ such that for all  $v^h \in S^h_0$,
\begin{equation}\label{fem2}
\int _{\Omega} \nabla u^h \cdot A\nabla v^h dx+\int _{\Omega} v^h b \cdot \nabla u^h dx+\int _{\Omega} c u^hv^h dx 
= \int_{\Omega} fv^h dx , 
\end{equation}
where $A=\left[ \frac{\partial \wG}{\partial r_{ij}}(w^h_k) \right]_{i,j=1}^2 $, 
$b=\left[ \frac{\partial \wG}{\partial p_{i}}(w^h_k) \right]_{i=1}^2 $, 
$c=\frac{\partial \wG}{\partial z}(w^h_k)$ and $f=G(u^h_k)$.

Under some additional assumptions on $G$, satisfied in particular by the Monge-Amp\`ere operator,
it is proved in \cite[Theorem 5.2]{BoehmerBook} and \cite[Theorem 9.1]{Boehmer08} 
that $u_k^h$ converges quadratically (as $k\to\infty$)  to a unique function $\hu^h$ satisfying the nonlinear equations
$$
( G(\hu^{h}),v^{h})_{L^2(\Omega)} =0 \quad \forall v^{h}\in S^h_0,$$
such that $\hu^h-u_0^h\in S^h_0$,
if the initial guess $u_0^h$ is close enough to $\hu^h$.
Moreover, $\hu^h$ converges to the solution $u$ of \eqref{e01} in $H^2$-norm as $h\to0$ if 
$u\in H^r(\Omega)$ for some $r>2$ and
the spaces $S^h_0$ possess appropriate approximation properties for functions vanishing on
$\partial\Omega$. Note that suitable approximation error bounds for the spaces $S^h_0=\sps$ 
have yet to be proved, see the results of \cite[Section 3]{DKS} for the spaces of continuous
piecewise polynomials vanishing on a piecewise conic boundary.
The stability of the MDS of Theorem~\ref{MDSCh5} and the dual local basis, related to the approximation power of the
space \cite{D07}, has been addressed in \cite{Abid_thesis}.

Note that in the case when $G$ is only conditionally
elliptic (e.g.\ elliptic only for a convex $u$ for Monge-Amp\`ere equation)
the ellipticity of the  linear problem \eqref{linear} is only guaranteed if $u^h_k$ 
satisfies the respective side condition $(x,u(x),\nabla u(x),\nabla^{2}u(x))\in \widetilde\Gamma$ 
for all $x\in\Omega$. For the Monge-Amp\`ere equation the side condition of convexity holds for
 $u^h_k$  if its second order derivatives are sufficiently close to those of the exact solution $\hu$.

\subsection{Implementation issues}

The standard techniques of the finite element method allow efficient computation of 
the solution $u^h$ of \eqref{fem2} using the local basis of $S^h_0=\sps$ described in 
Section~\ref{basis}. Moreover, efficient assembly algorithms for the polynomial Bernstein-B\'ezier shape
functions introduced in \cite{AAD11} can be employed in the same way as described in 
\cite[Section 5]{DKS} for the continuous polynomial finite elements on curved domains
enclosed by piecewise conics. We also refer to \cite{DSa13} for further implementation 
details related to fully nonlinear equations, and to \cite[Section 8]{ADS} for the efficient handling of the
global-local transformations in the finite element method relying on  Bernstein-B\'ezier shape functions.

\subsection{Numerical results}

In the experiments we focus on the Dirichlet problem for the prototypical and best studied 
Monge-Amp\`ere equation, 
\begin{equation}\label{MAch5}
G(u)=\det(\nabla^{2}u)-g= 0, \quad
 u|_{\partial\Omega} = \phi , 
\end{equation}
with $g(x)>0$, $x\in\Omega$, where the solution $u:\Omega\to\RR$ is assumed to be convex for 
the sake of uniqueness.  

We choose a number of test problems with a curved domain $\Omega$ bounded by piecewise conics,
a positive function $g$ and $\phi =0$. As in \cite[Section 6]{DKS}, starting from an initial triangulation of $\Omega$, we
obtain  a sequence of quasi-uniform triangulations $\tri^h$ by uniform refinement, whereby each triangle is subdivided into four
triangles by joining the midpoints of every edge. For each $h$, we use Böhmer's method described
above, with $S^h_0=\spsh$.	To solve \eqref{fem2} we use the 1-local basis corresponding to 
the MDS $M$ of Theorem~\ref{MDSCh5}.

We follow the suggestion of \cite[Remark 2.1]{DG06} to use an approximate 
solution of the Poisson problem
\begin{equation}\label{initial}
\Delta u = 2 \sqrt{g}, \quad
 u|_{\partial\Omega} = \phi , 
\end{equation}
as initial guess in the iterative schemes for the Monge-Amp\`ere equation \eqref{MAch5}. 
 Since $\phi =0$, we choose the initial guess $u^h_0$ in the same space $\spsh$ and obtain it by
 the standard Galerkin method.
 However, as  in \cite{DSa13}, we get much faster convergence of the Newton iteration 
 \eqref{newton} by a multilevel approach, where 
this initial guess is only used
on the initial triangulation, whereas on the refined triangulations 
a quasi-interpolant \cite[Section 5.7]{LSbook} of the last iterate from the 
previous level serves as an initial guess $u^h_0$. As a stopping criteria for Newton iterations \eqref{newton} on each level  the
following condition is employed:
\begin{equation}\label{term}
\|u^h_k -u^h _{k+1}\|_{L^2(\Omega)} < 10^{-15} .
\end{equation}

\begin{testproblem}\label{disk}\rm 
Equation \eqref{MAch5} in the unit disk $\Omega$ centered at the origin with  $g$ chosen 
such that the exact solution is $u=e^{0.5(x_1^2+x_2^2)}-e^{0.5}$.
\end{testproblem}

We use the same initial triangulation of the disk as in \cite[Example 2]{DKS}, see Figure~12 in
\cite{DKS}. 
The numerical results for Test Problem \ref{disk} are presented in Table~\ref{table1ch5},
which shows the $L^2$, $H^1$ and $H^2$ norms of the error $e_{\ell}=u_\ell-u$ of the 
last iterate $u_\ell=u^{h_\ell}_{m}$ on level $\ell$ against the exact solution and the
number $m$  of iterations \eqref{newton} for levels $\ell=1,\ldots,6$, where $\ell=1$ corresponds to
the initial triangulation. In addition, the first row of the table contains the errors of the
initial guess obtained by solving \eqref{initial} on the initial triangulation. The rate of
convergence between levels is estimated by the usual formula $\log_2(\|e_{\ell-1}\|/\|e_{\ell}\|)$.

The results show the convergence rates approaching  $O(h^6)$, $O(h^5)$ and $O(h^4)$ for the 
$L^2$, $H^1$ and $H^2$ norms, respectively, which is expectable since the solution $u$ is infinitely 
smooth and the space $\sps$ consists of piecewise polynomials of degree 5. The efficiency of the
multilevel approach to the computation of the initial guesses is also confirmed since  
only one or two Newton iterations are needed on each level to satisfy the termination criterion
\eqref{term}.

\begin{table}%
\begin{center} %
\begin{tabular}{|c|cc|cc|cc|c|}%
\hline
\hline %
$\ell$ &  $L^2$-error & rate  &  $H^1$-error & rate  &  $H^2$-error &  rate & 
$m$ \\ 
\hline %
init &  1.04e-2 &   &  3.20e-2 &  &  1.85e-1 & &  \\ 
\hline
1 &  2.12e-6 &  &  3.84e-5  &  &  1.25e-3 &  & 2 \\
2 &  2.98e-7 & 2.8 &  8.47e-6  & 2.2 &  3.35e-4 & 1.9 & 1  \\
3 &  6.79e-9 &  5.5  &  3.87e-7 & 4.5 & 2.86e-5 &  3.6 & 1  \\
4 & 1.36e-10 &  5.6 &  1.46e-8 & 4.7 &   2.12e-6 &  3.8 & 1  \\
5 &  2.52e-12  &  5.8  &  5.23e-10  &  4.8 &  1.47e-7 &  3.9 & 1  \\
6 &  9.51e-14 &  4.7  &  1.76e-11  &  4.9 &  9.53e-9 &  3.9 & 1 \\
\hline
\end{tabular}
\caption{Errors of the approximate solution and the rate of convergence for 
Test Problem \ref{disk} on the unit disk. $\ell$ indicates the level of refinement of the initial
triangulation, and $m$ is the number of Newton iterations \eqref{newton} on the $\ell$-th level.
The row marked `init' gives the errors of the initial guess on level 1.}
\label{table1ch5}
\end{center}
\end{table}

\begin{testproblem}\label{ellipse1}\rm 
Equation \eqref{MAch5} with $g(x)=e^{x_1}$ and $\phi=0$ in the elliptic disk  $\Omega$ with the boundary 
given by the equation $x_1^2+6.25x_2^2=1$. 
\end{testproblem}

The  initial triangulation is the same as the one used in \cite[Example 1 and Figure~7]{DKS}.
The results are presented in Table~\ref{table2ch5}. 
Since the exact solution $u$ is not known, we use alternative measures 
to estimate the error. One is the residual
\begin{equation}\label{residual}
R=\|G(u^{h}_k)\| _{L_2(\Omega)}, 
\end{equation}
and another is the  $L^2$, $H^1$ and $H^2$ norms of the difference 
$\varepsilon_\ell:= u_{\ell}-u_{\ell+1}$
between the approximate solutions $u_{\ell},u_{\ell+1}$ of two consecutive levels.
Note that in the case that $u_{\ell}$ converges to $u$ at least linearly in some norm, 
we may assume that
$\|u-u_{\ell+1}\| \leq \gamma\|u-u_{\ell}\|$
for some $\gamma < 1$ if $\ell$ is sufficiently large.
The triangular inequality then leads to 
$\|u-u_{\ell}\| \leq \frac{1}{1-\gamma}\|\varepsilon_\ell\|$, so that 
$\log_2(\|\varepsilon_{\ell-1}\|/\|\varepsilon_{\ell}\|)$ may serve as an estimate of the
convergence rate as long as it is positive.

We see that the numerical convergence rates in $L^2$, $H^1$ and $H^2$ norms are similar to those for Test
Problem~\ref{disk}. This indicates that the solution $u$ lies at least in $H^6(\Omega)$. In fact
it is expectable that $u$ should be infinitely differentiable because so are 
the data and the domain boundary. Note that \cite[Theorem~17.22]{GTbook} only assures that   
$u \in C^{2,\alpha}(\overline{\Omega})$ for all $0 < \alpha < 1$, but this
theorem only requires $C^3$ boundary and $C^2$ smoothness of $g$. The convergence rate of 
the residual \eqref{residual} is
close to $O(h^4)$, that is to the rate of the $H^2$-norm of the error, which is plausible because
$R$ is based on the second order derivatives of the approximate solution.

\begin{table}%
\begin{center} %
\begin{tabular}{|c|cc|cc|cc|cc|c|}%
\hline
\hline %
  $\ell$ & $\|\varepsilon_\ell\| _{L_2}$ & rate & 
 $\|\varepsilon_\ell\| _{H^1}$ & rate &  $\|\varepsilon_\ell\| _{H^2}$ & rate & $R$ & rate & $m$\\ 
\hline %
 init    &  &     &      &          &     &      & 6.58e-1 &  & \\ 
\hline %
 1   & 1.02e-8 &  & 3.64e-7 &  & 2.90e-5  &  & 4.95e-6 &  & 4 \\ 
2    & 9.59e-10 & 3.4 & 5.26e-8 & 2.8 & 6.37e-6  & 2.2 & 1.62e-6 & 1.6 & 1 \\ 
3    & 1.32e-11 & 6.2 & 1.29e-9 & 5.3 & 3.16e-7 & 4.3 & 1.37e-7 & 3.6 & 1 \\ 
4   & 2.25e-13 & 5.9 & 4.27e-11 & 4.9 & 2.05e-8 & 3.9 & 9.83e-9 & 3.8 & 1 \\ 
5    & 8.79e-15 & 4.7 & 1.61e-12 & 4.7 & 1.56e-9 & 3.7 & 6.61e-10 & 3.9 & 1 \\ 
6    & --- &     & --- &  & ---  &  & 4.33e-11 & 3.9 & 1  \\ 
\hline
\end{tabular}
\caption{Estimated errors of the approximate solution and the rate of convergence for 
Test Problem \ref{ellipse1} with $g(x)=e^{x_1}$ on the elliptic disk. The meaning of 
$\ell$, $m$ and `init' is the same as in Table~\ref{table1ch5}, 
$R$ is the residual error \eqref{residual} for the last iterate $u_\ell=u^h_m$ on level $\ell$, and
$\varepsilon_\ell:= u_{\ell}-u_{\ell+1}$ is the difference
between the approximate solutions of two consecutive levels. We left the entries
for $\ell=6$ related to $\varepsilon_\ell$ blank because their computation requires
the approximate solution $u_7$ of the next level.}
\label{table2ch5}
\end{center}
\end{table}

\begin{testproblem}\label{ellipse2}\rm 
Equation \eqref{MAch5} with $g(x)=\sin(\pi\vert x_1 \vert)+1.1$ and $\phi=0$ in the same 
elliptic disk $\Omega$ as in Test Problem~\ref{ellipse1}. 
\end{testproblem}

The numerical results can be found in Table~\ref{table3ch5}. 
Now \cite[Theorem~17.22]{GTbook} is not applicable because $g\notin C^2(\overline{\Omega})$.
Nevertheless, the method converges with approximate orders $O(h^{2.5})$, $O(h^{2.5})$ and 
$O(h^{1.5})$ for the $L^2$, $H^1$ and $H^2$ norms, respectively. This indicates that 
$u$ should be in $H^r(\Omega)$ for $r\approx 3.5$, but the approximation order of the method in 
$L^2$ norm is suboptimal.

\begin{table}%
\begin{center} %
\begin{tabular}{|c|cc|cc|cc|cc|c|}%
\hline
\hline %
  $\ell$ & $\|\varepsilon_\ell\| _{L_2}$ & rate & 
 $\|\varepsilon_\ell\| _{H^1}$ & rate &  $\|\varepsilon_\ell\| _{H^2}$ & rate & $R$ & rate & $m$\\ 
\hline %
 init    &   &  & & &  &  & 1.06e+0 & &  \\ 
\hline %
 1   & 2.92e-5 &  & 9.88e-4 &  & 9.48e-2 &  & 1.92e-2 & & 3 \\ 
2    & 5.41e-6 & 2.4 & 6.20e-5 & 3.9 & 4.44e-3  & 4.4 & 6.23e-3  & 1.6 & 2\\ 
3    & 1.21e-6 & 2.2 & 1.19e-5 & 2.4 & 1.40e-3  & 1.7 & 2.03e-3  & 1.6 & 1 \\ 
4     & 6.84e-8 & 4.1 & 2.01e-6 & 2.6 & 4.90e-4 & 1.5 & 7.46e-4 & 1.4 & 1 \\ 
5     & 1.44e-8 & 2.3 & 3.67e-7 & 2.5 & 1.47e-4 & 1.7 & 2.47e-4 & 1.6 & 1 \\ 
6     &  ---   &  &  ---   &  &  ---   &  & 9.04e-5 & 1.5 & 1 \\ 
\hline
\end{tabular}
\end{center}
\caption{Estimated errors of the approximate solution and the rate of convergence for 
Test Problem \ref{ellipse2} with $g(x)=\sin(\pi\vert x_1\vert)+1.1$ on the elliptic disk. 
The layout is the same as in Table~\ref{table2ch5}.}
\label{table3ch5}
\end{table}

\begin{testproblem}\label{C1domain}\rm 
Equation \eqref{MAch5} with $g(x)=1$ and $\phi=0$ in a $C^1$ domain
$\Omega$  bounded by the straight lines $x_2= \pm 1$ and semi-circles 
$$ 
x_1= \pm \left( 1+\sqrt{1-x_2^2}\right) ,\quad -1 \leq x_2 \leq 1.$$
\end{testproblem}

\begin{figure}[htbp!]
\centering
\begin{tabular}{c}
\includegraphics[height=0.45\textwidth]{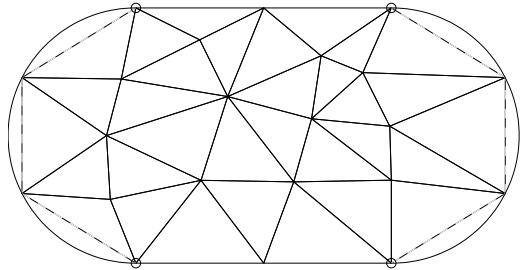} 
\end{tabular}
\caption{
The domain of Test Problem~\ref{C1domain} with initial triangulation. The boundary is $C^1$ 
at the four points marked with circles and $C^\infty$ elsewhere. 
Its top and bottom pieces are straight line segments.
}\label{InitMeshC1}
\end{figure}

The domain is visualized in Figure~\ref{InitMeshC1} together with the initial triangulation used
in our experiments. The straight line and circular segments are connected with $C^1$ smoothness
at the points $\pm (1,1)$ and $\pm (1,-1)$ indicated with circles. 

\begin{figure}[htbp!]
\centerline{
\begin{tabular}{c c}
\includegraphics[height=0.3\textwidth]{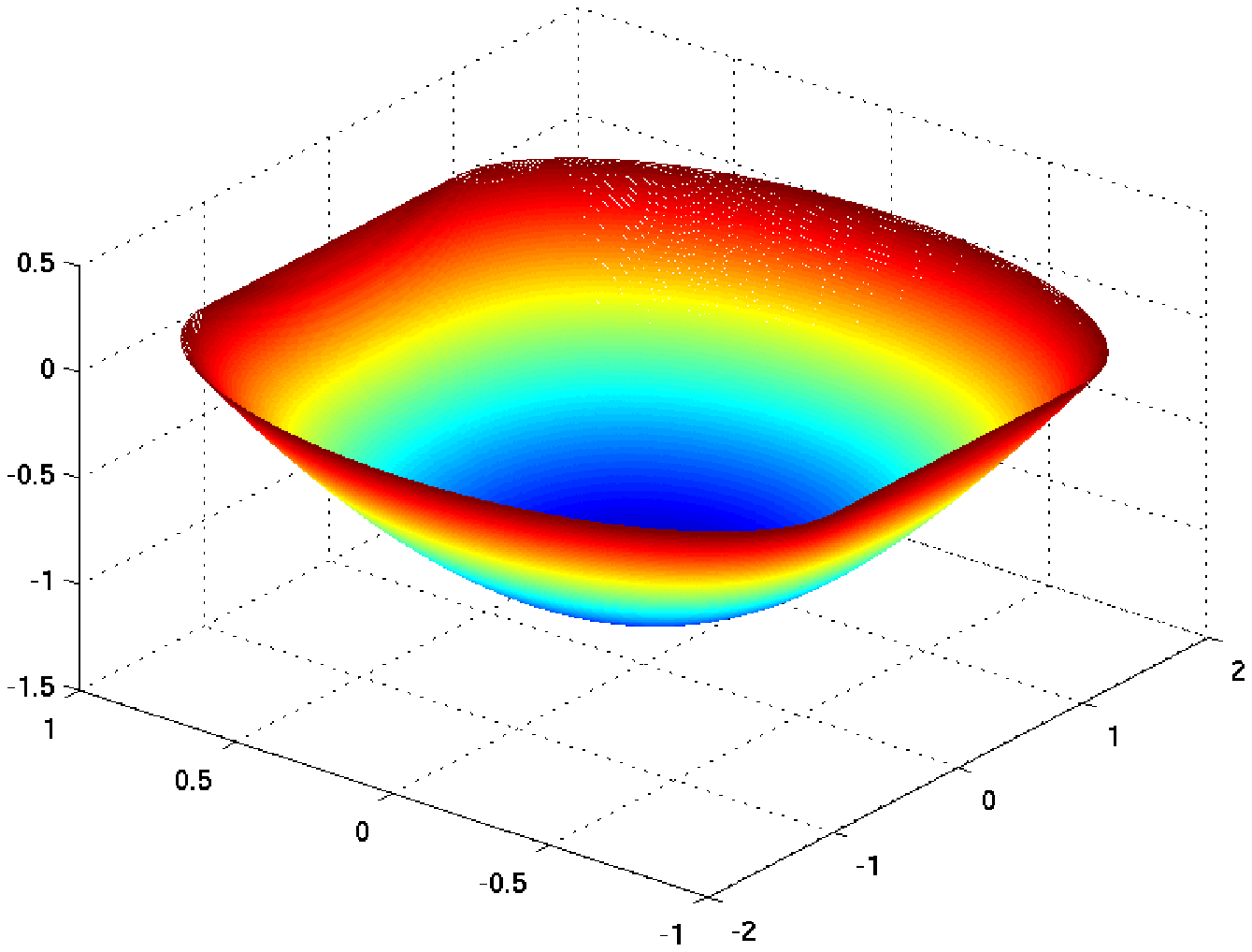} & \includegraphics[height=0.35\textwidth]{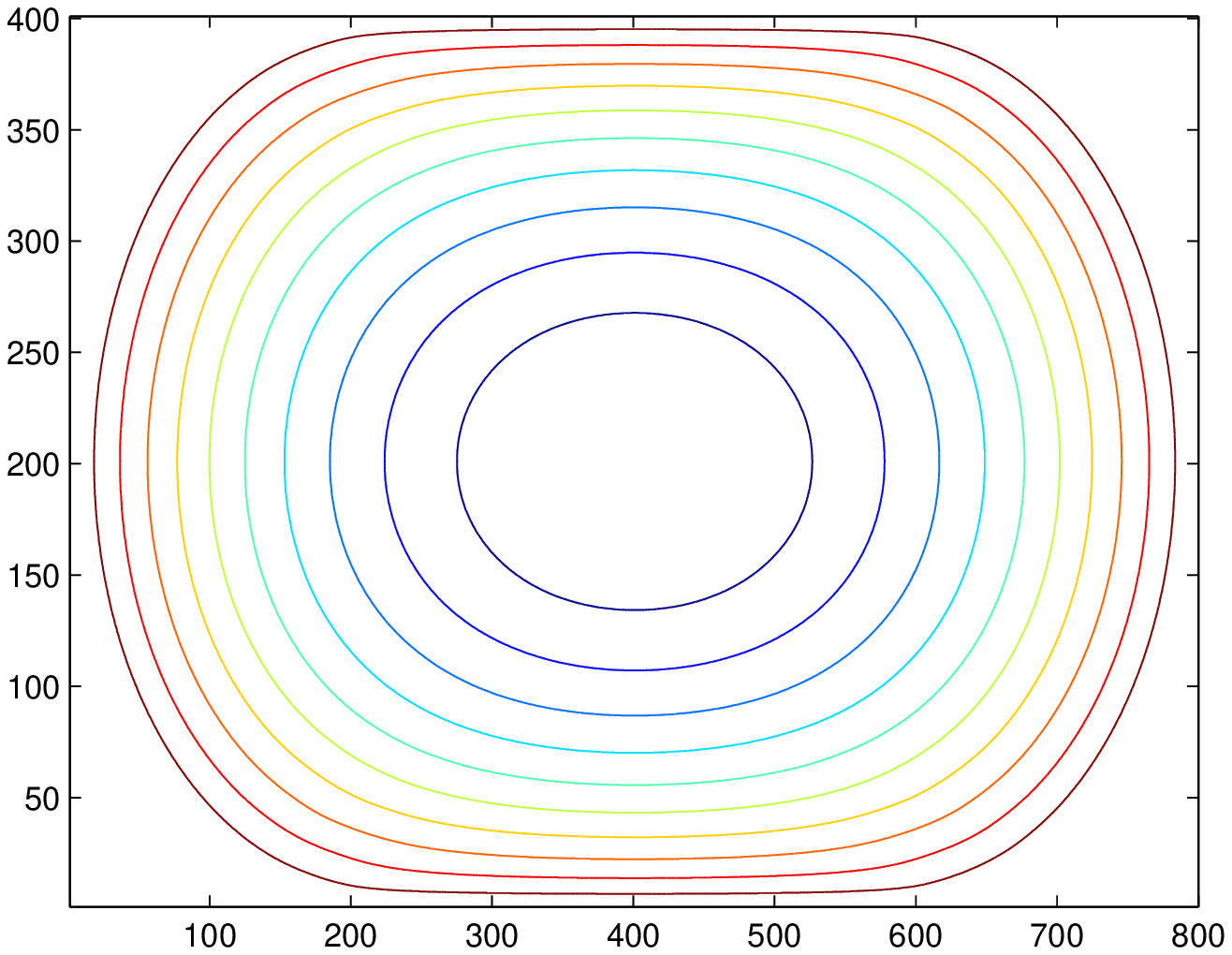}
\end{tabular}}
\caption{Approximate solution $u_\ell$ of Test Problem~\ref{C1domain} for the level $\ell=2$
 and its contor plot.}
\label{figtest3ch6}
\end{figure}

Similar to the tests with $g(x)=1$ on a square domain \cite[Section 5.1]{DSa13}, our
experiments do not show convergence of the method with respect to $\ell$. This is explained
in particular by the fact that the second derivatives of the solution $u$ of \eqref{MAch5}
with $\phi=0$ may not be continuous along any straight line boundary segment unless $g$
vanishes on this segment. Nevertheless, 
in contrast to the square domain, the approximate solutions $u_\ell$ keep the convex shape and 
the Newton iterations converge on each level. Figure~\ref{figtest3ch6} shows  $u_2$ and its
contor plot.

\begin{figure}[htbp!]
\centering
\begin{tabular}{c}
\includegraphics[height=0.45\textwidth]{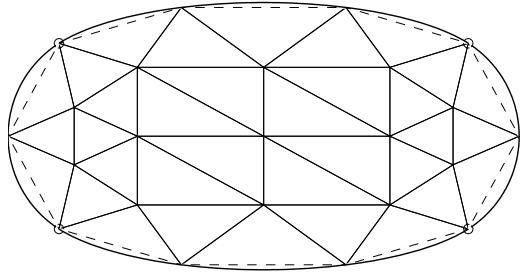} 
\end{tabular}
\caption{
The domain of Test Problem~\ref{C2domain} with initial triangulation. The boundary is $C^2$ 
at the four points marked with circles and $C^\infty$ elsewhere. 
}\label{InitMeshC2}
\end{figure}

\begin{testproblem}\label{C2domain}\rm 
Equation \eqref{MAch5} with $g(x)=1$ and $\phi=0$ in a centrally symmetric $C^2$ domain
$\Omega$  bounded by two elliptic and two circular segments,
see Figure~\ref{InitMeshC2}, where the top elliptic segment is given parametrically by the
equations
$$
x_1=4\cos t,\; x_2=1.3\sin t-c_2,\quad 0.15\pi\le t\le 0.85\pi,$$
and the left circular segment has radius $r$ and center $(c_1,0)$, with $r$ and $(c_1,c_2)$
being the radius and the center of the osculating circle to the ellipse 
$x_1=4\cos t$, $x_2=1.3\sin t$ at the point defined by $t=0.85\pi$.
\end{testproblem}

It is easy to check that elliptic and circular segments of $\Omega$ join with continuous 
curvature. We use the initial triangulation shown in Figure~\ref{InitMeshC2}.
The numerical results presented in Table~\ref{table4ch5} indicate 
$O(h^{4})$, $O(h^{3})$ 
and $O(h^2)$ convergence order in the $L_2$, $H^1$ and $H^2$-norm, respectively,
so that the solution $u$ is expected to belong to $H^r(\Omega)$ for $r\approx 4$. 
Note that \cite[Theorem 17.22]{GTbook} is not applicable because the boundary is not $C^3$.

\begin{table}%
\begin{center} %
\begin{tabular}{|c|cc|cc|cc|cc|c|}%
\hline
\hline %
  $\ell$ & $\|\varepsilon_\ell\| _{L_2}$ & rate & 
 $\|\varepsilon_\ell\| _{H^1}$ & rate &  $\|\varepsilon_\ell\| _{H^2}$ & rate & $R$ & rate & $m$\\ 
\hline %
 init    &   &  & & &  &  & 2.01e+0 & &  \\ 
\hline %
1   & 1.07e-3 &  & 1.00e-2 &  & 1.34e-1 &  & 9.10e-2 & & 2 \\ 
2    & 4.87e-5 & 4.5  & 8.56e-4 & 3.5  & 2.20e-2  & 2.6   & 2.20e-2  & 2.0 & 1 \\ 
3     & 3.04e-6 &  4.0  & 1.04e-4 & 3.0   & 5.30e-3 & 2.0  & 5.87e-3 & 1.9 & 1 \\ 
4     & 2.09e-7 & 3.7  & 1.39e-5 & 2.9   & 1.38e-3 & 1.9  & 1.56e-3 & 1.9 & 1 \\ 
5     & 1.58e-8  & 3.7 & 2.01e-6  & 2.8 & 3.80e-4  & 1.9 & 4.15e-4 & 1.9 & 1 \\ 
6     & ---  &  & ---  &  & ---  &  & 1.11e-4 & 1.9 & 1 \\ 
\hline
\end{tabular}
\end{center}
\caption{Estimated errors of the approximate solution and the rate of convergence for 
Test Problem \ref{C2domain}  on a $C^2$ domain. 
The layout is the same as in Table~\ref{table2ch5}.}
\label{table4ch5}
\end{table}


\begin{thebibliography}{44}

\def\JAT{J.~Approx.\ Theory}
\def\JCAM{J.~Comput.\ Appl.\ Math.}
\def\AiCM{Advances in Comp.\ Math.}
\def\CA{Constr.\ Approx.}
\def\SJNA{SIAM J. Numer.\ Anal.}
\def\CAGD{CAGD}

\bibitem{AAD11}
M. Ainsworth, G. Andriamaro and O. Davydov, Bernstein-B\'ezier  
finite elements of arbitrary order and optimal assembly procedures, 
{\sl SIAM J. Sci. Comp.}, {33} (2011), 3087--3109.

\bibitem{AAD15}
M. Ainsworth, G. Andriamaro and O. Davydov, A Bernstein-B\'ezier  
basis for arbitrary order Raviart-Thomas finite elements,
Constr. Approx. 41 (2015), 1--22.

\bibitem{ADS}
M. Ainsworth, O. Davydov and L. L. Schumaker, Bernstein-B\'ezier finite elements on
tetrahedral-hexahedral-pyramidal partitions, {\sl Computer Methods in Applied Mechanics 
and Engineering}, to appear. Preprint available from 
\url{https://www.staff.uni-giessen.de/odavydov/pyramids.html}

\bibitem{GA}
 G. Awanou, 
Pseudo transient continuation and time marching methods for Monge-Amp\`ere type equations,
 {\sl Advances in Computational Mathematics}, 41 (2015), 907--935. 
 
 \bibitem{MB93} 
M. Bernadou,
{Curved finite elements of class $C^1$: Implementation and numerical experiments. 
Part 1: Construction and numerical tests of the interpolation properties}, 
Comput. Method Appl. Mech. Engrg., 106(1-2) (1993), pp. 229--269.

\bibitem{ImplSurf} 
J. Bloomenthal et al, Introduction to Implicit Surfaces, Morgan-Kaufmann Publishers Inc., San Francisco, 1997.


\bibitem{Boehmer08} 
K. B\"ohmer, 
{On finite element methods for fully nonlinear elliptic equations of second order},  
{\sl SIAM J. Numer. Anal.}, 46(3) (2008), 1212--1249.

\bibitem{BoehmerBook} 
K. B\"ohmer, 
{Numerical Methods for Nonlinear Elliptic Differential Equations}: A Synopsis,
Oxford University Press, Oxford, 2010.

\bibitem{NBGS}
 S.C. Brenner, T. Gudi, M. Neilan, L.-Y. Sung, 
\emph{$C^0$ penalty methods for the fully nonlinear Monge-Amp\`ere equation}, Math. Comput., 80(276) (2011), 
1979--1995.

 
\bibitem{BrennerScott} 
S. C. Brenner,  and L.R. Scott,  {The Mathematical Theory
of Finite Element Methods,} Springer, New York, 1994.

\bibitem{D01}
O. Davydov,  Stable local bases for multivariate spline spaces,  {\sl  \JAT},
{111} (2001), 267--297.

\bibitem{D07}
O.~Davydov, {Smooth finite elements and stable splitting}, 
Berichte ``Reihe Mathematik'' der Philipps-Universit\"at Marburg, 2007-4 (2007). 
An adapted version has appeared as \cite[Section 4.2.6]{BoehmerBook}.


\bibitem{DKS}
O. Davydov, G. Kostin and A. Saeed, Polynomial finite element method for domains enclosed by piecewise
conics, {\sl  \CAGD}, to appear. 
\url{doi:10.1016/j.cagd.2015.11.002}
\url{arXiv:1510.00849}

\bibitem{DSa12}
O. Davydov and A. Saeed, Stable splitting of bivariate spline spaces
by Bernstein-B\'ezier methods, {\sl in}
``Curves and Surfaces - 
7th International Conference, Avignon, France, June 24-30, 2010'' (J.-D. Boissonnat et al, Eds.),
LNCS 6920,  Springer-Verlag, 2012,
pp.~220--235. %

\bibitem{DSa13}
O. Davydov and A. Saeed,  Numerical solution of fully nonlinear
elliptic equations by B\"ohmer's method, {\sl J.~Comput.\ Appl.\ Math.},  
{254} (2013), 43--54. %

\bibitem{DG06}
E. J. Dean and R. Glowinski, 
{Numerical methods for fully nonlinear elliptic equations of the Monge-Amp\`ere type}, 
{\sl Computer Methods in Applied Mechanics and Engineering}, 195 (2006), 1344--1386.

\bibitem{XF09-1}
X. Feng, M. Neilan,
\emph{Mixed finite element methods for the fully nonlinear Monge-Amp\`ere equation based on the vanishing moment 
method},
SIAM J. Numer. Anal., 47(2) (2009) 1226--1250.

 
\bibitem{GTbook}
D. Gilbarg and N. S. Trudinger, 
{Elliptic Partial Differential Equations of Second Order},  Springer-Verlag, 
 Berlin, 2001.

\bibitem{HCB05} 
T. J. R. Hughes, J. A. Cottrel, Y. Bazilevs,
{Isogeometric analysis: CAD, finite elements, NURBS, exact geometry and mesh refinement}, 
 Comput. Methods Appl. Mech. Engrg., 194(2005) 4135-4195.

\bibitem{LSbook} 
M. J. Lai and L. L. Schumaker, 
{Spline Functions on Triangulations},  Cambridge University Press, 2007.

\bibitem{LP13}
O. Lakkis, T. Pryer, 
A finite element method for nonlinear elliptic problems
{\sl SIAM Journal on Scientific Computing} 35 (2013), A2025--A2045.


\bibitem{Neilan14}
M. Neilan, Finite element methods for fully nonlinear second order PDEs based on a 
discrete Hessian with applications to the Monge-Ampère equation,
{\sl Journal of Computational and Applied Mathematics} 263 (2014), 351--369-



\bibitem{OA08}
A. Oberman, 
{Wide stencil finite difference schemes for the elliptic Monge-Amp\`ere equations and functions of 
the eigenvalues of the Hessian}, Discrete Contin. Dyn. Syst. Ser B 10(1) (2008), 221--238. 


\bibitem{Abid_thesis}
A. Saeed,
{Bivariate Piecewise Polynomials on Curved Domains, with Applications to Fully Nonlinear PDE's},
PhD thesis, University of  Strathclyde, Glasgow, 2012.

\bibitem{Sch15}
L. L. Schumaker,  Spline Functions: Computational Methods, SIAM (Philadelphia), 2015.

\end{thebibliography}
\end{document}